\documentclass[a4paper,twoside]{article}

\usepackage{amsmath, latexsym, amsfonts, amssymb, amsthm, amscd}
\usepackage{graphics,epsf,psfrag}
\usepackage[french,english]{babel}
\usepackage{titlesec}
\usepackage{lipsum}

\numberwithin{equation}{section}

\def \R {\mathbb{R}}
\def \N {\mathbb{N}}
\def \P {\mathbb{P}}
\def \E {\mathbb{E}}

\newtheorem{teo}{Theorem}[section]
\newtheorem{defi}[teo]{Definition}
\newtheorem{prop}[teo]{Proposition}
\newtheorem{lem}[teo]{Lemma}

\newtheorem{cor}[teo]{Corollary}

\newcounter{c}
\newcounter{aux}
\newcounter{aux1}
\newcounter{aux2}
\renewcommand{\thec}{\stepcounter{c}\arabic{c}}

\titleformat{\section}{\normalfont\scshape\center}{\thesection}{1em}{}
\titleformat{\subsection}{\normalfont\scshape}{\thesubsection}{1em}{}
\titleformat{\subsubsection}{\normalfont\scshape}{\thesubsubsection}{1em}{}

\begin{document}
\title{Front Velocity and Directed Polymers in Random Medium }

\author{Aser Cortines}
\date{}

\maketitle
\begin{abstract}
We consider a stochastic model of $N$ evolving particles studied by Brunet and Derrida \cite{Brunet2004, Cook1989}. This model can be seen as a directed polymer in random medium with $N$ sites in the transverse direction. In \cite{Cook1989}, Cook and Derrida, use heuristic arguments to obtain a formula for the ground state energy of the polymer. In this paper, we formalize their argument and show that there is an additional term in the formula in the critical case. We also consider a generalization of the model, and show that in the noncritical case the behavior is basically the same, whereas in the critical case a new correction appears. 
\end{abstract}

\section{Introduction}\label{sec:model.def}

A relatively simple formulation for the problem of directed polymer in random medium is the following \cite{Cook1989, Giacomin2007}. The lattice consists of $L$ planes in the transversal direction. In every plane there are $N$ points that are connected to all points of the previous plane and next one. For each edge $ij\,$, connecting the $t$-th plane to the $(t+1)$-th plane, a random energy $\xi_{ij} (t+1)$ is sampled from a common probability distribution $\xi$. With a slight abuse of notations we write $\xi$ for both the distribution and a random variable with distribution $\xi$. For $\omega = [\,\omega_1 \,,\ldots ,\omega_L \,]$ a standard random walk on $\mathbb{G}_N$ the complete graph on $N$ vertices, we define the energy $E_{\omega}\,$ of the directed path by summing the energies of the visited bounds
$$
E_{\omega} : = \sum_{s=1 }^{L} \xi_{\omega_s \, \omega_{s+1} } (s+1) \,.
$$
We define the probability measure $\mu_L$ on the space of all directed paths of length $L$ by 
$$
\mu_L (\omega ):= Z_L (T)^{-1} \exp ( -E_{\omega} / T ) \, ,
$$
where $T$ is the temperature and $Z_L (T)$ is the partition function. The directed path $\big(\omega_i , i \big)_{i \geq 0}$ can be interpreted as a polymer chain living on $\mathbb{G}_N\times \N$, constrained to stretch in one direction and governed by the Hamiltonian $\exp ( - E_{\omega} / T)$\,.

We will be interested in the case where the random energies $\xi$ depends on $N$ the number of vertices of $\mathbb{G}_N$ and we will work at zero temperature. When $T=0$, we are faced with an optimization problem: computing the ground state energy of the model  {\it i.e.} the lowest energy of all possible walks.

In \cite{Cook1989} Cook and Derrida consider the particular case of zero temperature and $\xi$ distributed according to a Bernoulli of parameter $1/N^{1+r}$, with $r\geq 0$, which they call the percolation distribution. Hence, the energy $E_\omega$ of a directed path of length $L$ is equal to the number of times $\xi_{\omega_s \, \omega_{s+1} } = 1$ along this path. Moreover, we can easily conclude that if the ground state of the polymers of length $L$ is $E_L$, then $E_{L+1} \leq E_{L}+1$. 

For $N$ fixed, the ratio $E_L / L$ converges and is a constant a.s. In \cite{Cook1989} the authors call this limit the ground state energy per unity of length and they derive the following asymptotic for it, when $N \to \infty$ 
\begin{equation} \label{equa:bad.formula}
E = \Big(1+\lfloor 1/r \rfloor \Big)^{-1},
\end{equation}
where $\lfloor \cdot \rfloor$ denotes the integer part. Their statement is based on the observation that the typical number of sites on the $t$-th plane connected to the first plane by a path of zero energy  is $N^{1-tr}$. Hence, if $N$ is large enough and $1-tr$ positive there is a path of zero energy (which is necessarily a ground state) from $0$ to $t$, whereas when $1-tr$ is negative there is no such path. Their argument, although informal, is correct, but the case where $1/r$ is an integer (the critical case) requires a more careful analysis. In this paper we formalize their argument and show that there is an additional term in (\ref{equa:bad.formula}) when $1/r$ is an integer.
\medskip

In this paper, we choose to approach the polymer problem described above through the point of view of an interacting particles system. It consists in a constant number $N$ of evolving particles on the real line initially at the positions $X_1(0), \ldots ,X_N(0)$. Then, given the positions $X_i(t)$ of the $N$ particles at time $t \in \N$, we define the positions at time $t + 1$ by:
\begin{equation} \label{definition.X.derrida.brunet}
X_i (t + 1) : = \max_{1\leq j \leq N} \big\{ X_j(t) + \xi_{j,i} (t + 1) \big\} ,
\end{equation}	
where $\big\{ \xi_{i,j} (s) \, ; 1 \leq i, j \leq N\,, s \in \N \big\} $ are i.i.d. real random variables of common law $\xi$. The $N$ particles can also be seen as the fitness of a population under reproduction, mutation and selection keeping the population size constant. 

Moving fronts are used to model some problems in biology and physics. It describes, for example, how the fitness of gene propagates through a population. In physics they appear in non-equilibrium statistical mechanics and in the theory of disordered systems \cite{Derrida1988}.

We now explain how (\ref{definition.X.derrida.brunet}) is related to the polymer problem studied by Cook and Derrida in \cite{Cook1989}. One can check by induction that
$$
X_i(t)= \max \Big\{ X_{j_0}(0) + \sum_{s=1}^{t} \xi_{j_{s-1} j_{s}} (s); 1 \leq j_s \leq N, \ \forall s = 0, \ldots, t-1 \text{ and } j_t = i \Big\} \,.
$$
Then we take $X_j(0) = 0$ for all $1 \leq j \leq N$ and sample $-\xi_{ij}(t)$ as a Bernoulli of parameter $1/N^{1+r}$. From the above formula, $-X_i (t)$ corresponds to the ground state energy of the polymer conditioned to be on $i$ at $t$. Therefore, the ground state is obtained by taking the maximum over all possible positions.  

\begin{defi}[Front Speed] Let $\phi \big(X(t)\big) = \displaystyle \max_{1 \leq i \leq N} \ \big\{ X_i(t) \big\}$. The front speed $v_N$ is defined as 
\begin{equation}\label{equa.def.speed}
v_N : =\lim_{t \to \infty} \frac{\phi \big( X(t) \big)}{t} .
\end{equation}
For $N$ fixed, the limit (\ref{equa.def.speed}) exists and is constant a.s., see \cite{Comets2013} for more details and a rigorous proof.
\end{defi}
Hence, it is not difficult to see that the ground state energy per unit of length is equal to $-v_N$ as defined in (\ref{equa.def.speed}).

\medskip

This model was introduced by Brunet and Derrida in \cite{Brunet2004} to better understand the behavior of some noisy traveling-wave equations, that arise from microscopic stochastic models. By the selection mechanism, the particles remain grouped, they are essentially pulled by the leading ones, and the global motion is similar to a front propagation in reaction-diffusion equations with traveling waves. In \cite{Brunet2004}, Brunet and Derrida solve for a specific choice of the disorder ($\xi_{ij}$ are sampled from a Gumbel distribution) the microscopic dynamics and calculate exactly the velocity and diffusion constant. Comets, Quastel and Ramirez in \cite{Comets2013} prove that if $\xi$ is a small perturbation of the Gumbel distribution the expression in \cite{Brunet2006} for the velocity of the front remains sharp and that the empirical distribution function of particles converges to the Gumbel distribution as $N \to \infty$. They also study the case of bounded jumps, for which a completely different behavior is found and finite-size corrections are extremely small. 

Traveling fronts pulled by the farmost particles are of physical interest and not so well understood, see \cite{Panja2004} for a survey from a physical perspective. It is conjectured that, for a large class of such models where the front is pulled by the farmost particles, the motion and the particle structure have universal features, depending mainly on the tails distribution \cite{Brunet2004, Brunet2006}. Recent results have been rigorously proved for different models in front propagation confirming some of the conjectures. 

B\'erard and Gou\'er\'e  \cite{Berard2010} consider the binary Branching Random Walk (BRW) under the effect of a selection (keeping the $N$ right-most particles). They show that under some conditions on the tail distribution of the random walk the asymptotic velocity converges at the unexpectedly slow rate $(\log N)^{-2} \,$ . Couronn\'e and Gerin  \cite{Couronne2011} study a particular case of BRW with selection where the corrections to the speed are extremely small. Maillard in \cite{Maillard2011} shows that there exists a killing barrier for the branching Brownian motion such that the population size stays almost constant. He also proves that the recentered position of this barrier converges to a Levy process as $N$ diverges. In the case where there are infinitely many competitors evolving on the line, called the Indy-500 model, quasi-stationary probability measures are superposition of Poisson point processes \cite{Ruzmaikina2005}. 

\bigskip

In the first part of this paper, we study the model presented in \cite{Cook1989} and described above. We consider the case where the distribution of the $\xi_{ij}$ depends on $N$ and is given by
\begin{align}\label{equa:defi.2.states}
\P &\big( \xi (N) =0 \big)  = p_0(N)  \sim  \rho/N^{1+r}  \\
\P & \big( \xi (N) = - 1 \big) = 1 -  \P\big( \xi (N) =0 \big) \,, \nonumber
\end{align}
where $r>0$, $\rho>0$ and for sequences $a_N , b_N$ we write $a_N \sim b_N$ if $a_N / b_N \to 1\,$. We will often omit $N$ in the notation. Since $\xi$ is non-positive, the front moves backwards. As a consequence of the selection mechanism and the features of $\xi$, all particles stay from a distance at most one from the leaders. And when the front moves, {\it i.e.} $\phi\big(X(t)\big) = \phi\big(X(t-1)\big)-1$, all particles are at the same position. This particular behavior hides a renewal structure that will be used when computing the front speed.

The case $1/r \in \N$  is critical and the system displays a different behavior. For $N$ large enough, at time $t= 1/r$, we show that there is a Poissonian number of particles $X_i $ that remain in zero. Then, at the $1/r$-th plane there exists a finite number (possibly zero) of sites that can still be connected to the first plane through a path of zero energy. Whereas, when $1/r \not\in \N$ the typical number of such sites is of order $N^{1-tr}$. This difference of behavior leads to an additional term in (\ref{equa:bad.formula}) and the following Theorem holds.

\begin{teo}\label{teo:derrida} Let $\xi$ be distributed according to (\ref{equa:defi.2.states}). Then the front speed $v_N$ satisfies
\begin{equation}
\lim_{N \to \infty} v_N = \left\{ 
\begin{array}{lcl}
-\big( \, 1+ \lfloor   1/r   \rfloor \big)^{-1}, & \text{if} &  1/r  \not\in \N \\
-\big( \, 1+ \lfloor   1/r \rfloor - e^{- \rho^{1/r}} \,\big)^{-1} , & \text{if} &  1/r  \in \N \, ,
\end{array}
\right.
\end{equation}
In the case where $r=0$,
\begin{equation}
\lim_{N \to \infty} v_N = 0\,.
\end{equation}
\end{teo} 

In Section \ref{sec.front.speed.three.states} we consider the case where $\xi$ takes values in the lattice $\mathbb{Z}_0 = \{ l \in \mathbb{Z} ; l \leq 0 \}\,.$ Then we set for $i \in \N$
\begin{equation}\label{equa:defi.3.states}
p_i(N) = \P(\xi(N) =-i)  \,,
\end{equation} 
and assume that $p_0 \sim \rho / N^{1+r}$ where $r$ and $\rho$ are non-negative. Let
\begin{equation}\label{equa.defi.q.2}
q_2(N) : = \P\big( \xi (N) \leq -2 \big)=1-p_0-p_1\,.
\end{equation}
We also assume that for $i \geq 2 $
\begin{equation}\label{equa.defi.vartheta}
\frac{p_i(N)}{q_2(N)} = \P (\vartheta = -i)\,,
\end{equation}
where $ \vartheta $ is an integrable distribution on the lattice $\mathbb{Z}_{-2}$ that does not depend on $N$. 

\medskip

As we explain in Section \ref{sec.conclusion.teo.speed.3.states.gen}, we can further generalize the model and consider $\xi$ distributed as 
\begin{equation}\label{equa:defi.3.states.intro.0}
\xi = p_0(N) \delta_{\lambda_0} + p_{1}(N) \delta_{\lambda_1} +q_2(N) \vartheta   (dx) \,,
\end{equation} 
where $\vartheta (dx)$ is an integrable probability distribution over $]-\infty\,, \lambda_1 [$ and $\delta_{\lambda_i}$ is the mass distribution. We assume that $\lambda_1 < \lambda_0$ and that $p_0(N) \sim  \rho/N^{1+r}$ for a $r>0\,.$ Then, the velocity $v_N$ obeys the following asymptotic.


\begin{teo} \label{teo.speed.3.states.gen} Let $\xi$ be distributed according to (\ref{equa:defi.3.states.intro.0}). Assume that 
$$
p_0(N) \sim \frac{\rho}{N^{1+r}}, \qquad \text{and} \quad \lim_{N \to \infty} q_2(N) = \theta \,,
$$
where $r>0$ and $ 0<\theta<1$. Then the front speed $v_N$ satisfies
\begin{equation}
\lim_{N \to \infty} v_N = \left\{ 
\begin{array}{lcl}
\lambda_0- (\lambda_0-\lambda_1)\big( \,1+ \lfloor  1/r  \rfloor  \big)^{-1}, & \mbox{if} & 1/r \not\in \N   \\
\lambda_0 - (\lambda_0-\lambda_1)\big( \, \lfloor   1/r  \rfloor  + 1- 1/g(\theta) \, \big)^{-1} , & \mbox{if} & 1/r  \in \N \,,
\end{array}
\right.
\end{equation}
where $g (\theta ) \geq 1 $ is a non-increasing function. The conclusion in the case $1/r \not\in \N$ still holds if $\xi$ satisfies the weaker assumption $q_2/(  1-p_0  ) \leq \theta' \,,$ for some $0< \theta'<1\,$.
\end{teo}
\medskip

The paper is organized as follows.  In Subsection \ref{sec:number.particles.zero} we compute the typical number of leading particles, which corresponds to the number of paths of zero energy, and in Subsection \ref{subsec.front.speed}, we calculate the limit of $v_N$ as $N\to \infty$, exhibiting in particular the additional term appearing in (\ref{equa:bad.formula}) in the critical case. In Subsection \ref{sec:jumping} we compute the typical number of leading particles, when $\xi$ is distributed according to (\ref{equa:defi.3.states.intro.0}). Subsections \ref{bound} and \ref{sec.conv.integral} present some technical results and calculations. In Subsection \ref{sec:front.speed.3.states} we compute the front velocity and prove the discrete version of Theorem \ref{teo.speed.3.states.gen}. Finally, in Section \ref{sec.conclusion.teo.speed.3.states.gen} we sketch the proof of Theorem \ref{teo.speed.3.states.gen} .

\section{Front speed for the two-state percolation distribution}\label{sec.derrida.two}

As in \cite{Comets2013}, we consider the following stochastic process.

\begin{defi} Let $Z(t) := \big( Z_l (t) ; \,  l= 0,1 \big)$ be defined as
\begin{equation}
Z_l( t  ) = \sharp \big\{  i  ; 1 \leq j \leq N  ;  X_i( t ) = \phi ( X(t - 1)  ) - l  \big\} \, ,
\end{equation}
where $\sharp $ denotes the number of elements of a set. Recall that $\phi ( X(t - 1)  ) $ is the position of the front at $t-1$.
\end{defi}

Note that, due to the special features of the distribution (\ref{equa:defi.2.states}), $Z_0(  t )$ is equal to the number of leaders if the front has not moved backwards between times $t-1$ and $t$, and to $0$ if the front moved. $Z$ is a homogeneous Markov chain on the set
$$
\Omega(N) = \bigr\{ \, x \in \{  0, 1 ,  \ldots , N   \}^2  \, ; x_0 + x_{1} = N \, \bigl\} \, ,
$$ 
where $x_i$ are the coordinates of $x$. It is obvious that $Z_0$ determines completely the vector $Z (t)\,$, since $Z_1(t) = N - Z_0( t )$. The transition rates of the Markov chain $Z_0(t)$ is given by the Binomial distributions
\begin{align}\label{equa:law.Z}
\P & \big( Z_0 ( t+1  ) = \cdot  \mid Z(t) = x \big) \nonumber \\ 
& = \P \big( Z_0 (  t+1  ) = \cdot \mid Z_0( t  ) = x_{0} \big) =\left\{ 
\begin{array}{lc}
\mathcal{B} \big( N, 1-(1-p_0)^{x_{0}} \big)( \cdot), &  x_0 \geq 1 \\
\mathcal{B} \big( N, 1-(1-p_0)^{N} \big)( \cdot) , &  x_0 = 0 \,.
\end{array}
\right. 
\end{align}

We will often consider Markov chains with different starting distributions. For this purpose we introduce the notation $\P_\mu $ and $\E_\mu$ for probabilities and expectations given that the Markov chain initial position has distribution given by $\mu$. Often, the initial distribution will be concentrated at a single state $x$. We will then simply write $\P_x$ and $\E_x$ for $\P_{\delta_x}$ and $\E_{\delta_x}$. 

In this Section, $\oplus $ denotes the configuration $(0,N) \in \Omega(N)\,$. Furthermore, we introduce the notation
\begin{equation}\label{eq:integer.part.r}
1/r = m + \eta \,,
\end{equation} 
where $m$ stands for the integer part of $1/r$ and $\eta$ its fractional part.

\subsection{Number of Leading Particles} \label{sec:number.particles.zero}

In this Subsection, we show that under a suitable normalization and initial conditions the process $Z_0$ converges as $N$ goes to infinity. 

Consider the random variable 
\begin{equation} \label{eq:moving.backwards.stopping.time}
\tau = \inf \big\{t\,;\phi\big(X(t)\big) < \phi\big(X(t-1)\big) \,\big\} \,.
\end{equation}
Then, $\tau$ is a stopping time for the filtration $\mathcal{F}_t = \{ \xi_{ij} (s); s \leq t \,\,\, \text{and} \,\,\, 1\leq i,j\leq N \}$. It is not difficult to see that $\tau$ is also the first time when $Z_0$ visits zero
$$
\tau =  \inf \big\{t\,; Z_0 (t) = 0 \big\} \,,
$$
and that $Z$ starts afresh from $\oplus\, $ when $t = \tau $ {\it i.e.} the distribution of $Z(\tau +t)$ is the same as the distribution of $Z(t)$ under $\P_\oplus$. 

\begin{defi} \label{defi:Y} 
Let $Y( t )$ be the number of leading particles at time $t$ if the front has not moved
\begin{equation}
Y( t ) := Z_{0} (t) 1_{\{t \leq \tau \}} \,.
\end{equation}
\end{defi}
Then, $Y$ is a homogeneous Markov chain with absorption at zero and transition rates given by the Binomial distributions
$$
 \P \big(Y( t+1 ) = \cdot   \mid Y( t ) = k \big)= \mathcal{B} \big( \, N,  1-( 1-p_0)^k \, \big) (\cdot) \,.
$$
The advantage of working with $Y$ rather than $Z_0$ is that the above formula holds even if $Y(t) = 0$. 

\begin{prop}\label{prop:gen.function} Let $\xi$ be distributed according to (\ref{equa:defi.2.states}). For $k \in\{1,2,\ldots N \}$ denote by $G_{k}( s,t )$ the Laplace transform of $Y(t)$ under $\P_{k}$ at $s\in \R$. Then,
\begin{equation}\label{equa:prop.gen.function}
G_{k}(s,t): = \E_{k} [e^{s \, Y(t)}]= \exp \Big\{( e^s-1 )k \big(N p_0\big)^t \big(1 + o\big( 1\big) \big) \Big\} \, .
\end{equation}
\end{prop}
\emph{Proof.} Conditioning on $\mathcal{F}_{t-1}: = \{\xi_{ij} (s); s\leq t-1\} \,,$ 
\begin{align*}
\E_{k} [ e^{s \, Y( t )} ] & =\E_{k} \big[ \, \E[e^{s\,Y( t)} \mid Y(t-1)] \, \big] \\
& = \E_{k} \Big[ \Big( 1 + (  e^s-1 ) \big( 1- ( 1- p_0 )^{Y( t-1)} \big)\Big)^N \,\Big]\,.
\end{align*}
Since $p_0 \sim \rho/N^{1+r}$ with $r>0$ and $Y(t-1) \leq N $ we obtain by first order expansion that 
$$
\Big( 1  + (e^s-1) \,\big(1- (1-p_0 )^{Y(t-1)}\big) \Big)^N = s_{(1)}(N)^{Y(t-1)} \, ,
$$
where $\displaystyle s_{(1)}(N)= \exp \big\{ (e^s - 1) \big( N p_0 + o (N p_0) \big) \big\} \,  $ and $o ( N p_0  )$ converges to $0$ independently from $Y(t-1)$. 

Repeating the argument,
$$
\E_{k} [ \ e^{s\,Y(t)} \ ]=  \E_{k} [s_{(1)}(N)^{Y(t-1)}] = \E_{k} [s_{(2)}(N)^{Y( t-2)}] \, ,
$$
\vspace{0.2cm}
with $\displaystyle s_{(2)}(N)= \exp \big\{ (s_{(1)}(N)-1) \big( N p_0 + o (  N p_0 )\big) \big\} \,  $. 

Expanding $s_1(N)-1$,
\vspace{0.2cm}
\begin{align*}
s_{(1)}(N)-1 & = \exp \Big\{ (e^s - 1) \big( N\, p_0 + o (   N p_0  ) \big) \Big\} - 1 \\
& = (e^s-1) \big( N\, p_0 + o (N p_0)  \big) \, .
\end{align*}
Hence, $s_{(2)}(N)= \exp\big\{(e^s-1) \big(N\, p_0\big)^2 + o \big( (N p_0)^2\big) \big\} \,$. We proceed recursively and obtain the expression
$$
\E_{k} [e^{s\,Y(t)}]= \exp  \Big\{ k \,( e^s-1) \big(N\, p_0 \big)^t \big(  1 + o (1)  \big) \Big\} \,,
$$ 
which proves the statement. \hfill $\Box$
\medskip
 
We point out that the case $k = N \, $ corresponds to $Z(0) = \oplus \, $. We now state  two Corollaries of Equation (\ref{equa:prop.gen.function}).

\begin{cor} \label{cor:Y.big.t} Let $\xi$ be distributed according to (\ref{equa:defi.2.states}) and $k \in \{1, \ldots, N\}$. Then, for $t \geq m+1 $ ,
\begin{equation}\label{equa.cor:Y.big.t}
\P_{k} \big(Y(t) = 0\big) \geq  1 - \rho^{t} \, N^{1-t\,r} + o\big(N^{1-t\,r}\big) \, .
\end{equation}
\end{cor}
\emph{Proof.} Since $\P_{k} \big( Y(t) = 0\big) = \lim_{s \to - \infty } E_{k} \big[ e^{s Y(t)} \big]$, Proposition \ref{prop:gen.function} implies that
$$
\P_{k} \big( Y(t) = 0\big) = \exp \Big\{- k \big(N p_0\big)^t \big(1 + o( 1) \big) \Big\} \geq \exp \Big\{ - N \big(N p_0\big)^t \big(1 + o( 1) \big) \Big\} \,.
$$
Then, we obtain (\ref{equa.cor:Y.big.t}) by first order expansion. \hfill $\Box$


\vspace{0.2cm}

\begin{cor} \label{cor:converg.Y.eta.random} Let $\xi$ be distributed according to (\ref{equa:defi.2.states}) with $\eta=0$, {\it i.e.} $r=1/m \,$. Assume that $\kappa (N)$ is a sequence of random variables in $\,\{1, \ldots, N\}\,$ that are independent from $\xi_{ij}$ and that $\kappa(N)/N$ converges a.s. to a positive random variable $U$.  

Then, under $\P_{\kappa(N)}$, $Y(m)$  converges in distribution to $ Y_\infty $ a doubly stochastic Poisson random variable characterized by its Laplace transform
\begin{equation}
\E[ e^{s\,Y_\infty }]:= \E \big[  \exp \big\{ U \,(e^s-1)\rho^{m} \big \} \big].
\end{equation}
\end{cor}
\emph{Proof.} From (\ref{equa:prop.gen.function}), we have that
$$
\E [ e^{s Y(m)} \mid \kappa (N)  ]  = \exp \big\{ ( e^s-1 ) \, \kappa(N) \rho^{m} \, N^{-1} \big(1 + o(1) \big) \big\} \,.
$$ 
The term $o(1)$ converges to zero independently from $\kappa(N)\,$ the initial position. Then, by dominated convergence, we obtain that
$$
\lim_{N \to \infty }\E [ e^{s Y(m)}  ]  = \E \big[  \exp \big\{ U \,(e^s-1)\rho^{m} \big\} \big] \,,
$$
which concludes the proof. \hfill $\Box$
\vspace{0.2cm}

We now prove a large deviation principle for $Y$. As in \cite{Dembo2010, DenHollender2008}, we denote by
\begin{equation}
\Lambda_{k,\,t} (s) := \lim_{N \to \infty} \frac{1}{k N^{-rt}} \log \E_{k } \big[ e^{s \, Y(t)} \big] \,,
\end{equation}
the cumulant generating function of $Y$ under $\P_k$. From (\ref{equa:prop.gen.function}) we see that $\Lambda_{k,\,t} (s) = (e^s-1)\rho^{t}\, $. Denoting by
\begin{equation}
\Lambda_{k,t}^* (x) := \sup_{s \in \R} \{ x s - \Lambda_{k, t} (s)\}\,,
\end{equation}
the Legendre transform of $Y(t)$ under $\P_{k}$, we have that
\begin{equation}\label{equa:rate.function.Y}
\Lambda_{k,\,t}^* (x) = \left\{ 
\begin{array}{ccl}
x ( \log x -\log \rho^t ) +\rho^t  -x \,, & \mbox{ if } & x > 0 \\ 
\infty \,, & \mbox{ if } & x \leq 0 \,.\\ 
\end{array}
\right.
\end{equation}
\vspace{0.2cm}

\begin{prop}(Large Deviation Principle for $Y\,$)\label{prop:large.deviation.Y}  Let $\xi$ be distributed according to (\ref{equa:defi.2.states}). For $t \leq m$, let $k(N) \leq N $ be a sequence of positive integers such that
$$
\lim k(N) \, N^{-r\,t} = \infty \,.
$$
Then, under $\P_{k(N)}$, $Y(t)/\big(k(N)N^{-r\,t}\big)$  satisfies a Large Deviation Principle with rate function given by $\Lambda^*_{k,t}$ as in (\ref{equa:rate.function.Y}) and speed $k(N)N^{r\,t}$.
\end{prop}

\emph{Proof.} In fact, it is a direct application of G\"artner-Ellis Theorem (see {\it e.g.} Theorem V.6 in \cite{DenHollender2008}). Since $\Lambda$ is smooth, it is a lower semi-continuous function, therefore the lower bound in the infimum can be taken over all points. \hfill  $\Box$
\vspace{0.2cm}

Our next Corollary formalizes the statement of Cook and Derrida in \cite{Cook1989}. 



\begin{cor} \label{cor:fluct.Y.random}  Let $\xi$ be distributed according to (\ref{equa:defi.2.states}) with $\eta>0$. Assume that $\kappa(N)$ is a sequence of random variables in $\,\{1, 2, \ldots, N\}\,$ independent from $\xi_{ij}$ and that $\kappa(N)/N$ converges a.s. to $U$ a positive random variable. 

Then, under $\P_{\kappa(N)}$ for $t \leq m$
\begin{equation}
\lim_{N\to \infty } \P \Big( \ \Big| \,\frac{ Y(t) }{ \rho^t\, U\, N^{1-t r} } -1 \,  \Big| \geq \varepsilon \ \Big) =  0 \, .
\end{equation}
The same statement holds for $\eta = 0$ and $t \leq m-1 $.
\end{cor}

\emph{Proof.} We first consider the case where $\kappa(N)$ is a deterministic sequence and $\kappa(N)/N \to u \,,$ with $0 < u \leq 1$. Then the conditions of Proposition \ref{prop:large.deviation.Y} are satisfied and $Y(t)/\big( \kappa(N)\,N^{-t\,r} \big) $ satisfies a Large Deviation Principle with rate function given by (\ref{equa:rate.function.Y}), which only zero is at $\rho^t$. This implies the desired convergence.  

The random case is solved by conditioning on $\kappa(N) = Y(0)$.
\begin{align*}
\P & \Big( \ \Big|\, \frac{ Y(t) }{ \rho^t U N^{1- t r} } -1 \, \Big| \geq \varepsilon \ \Big) \\
 & = \int \P^{(2)}_{\kappa(N)(\omega_1)}\Big( \, \Big| \, \frac{ Y(t) }{ \rho^t U(\omega_1) N^{1- t r} } -1 \, \Big| \geq \varepsilon \, \Big) \P^{(1)} (\mathrm{d} \omega_1) \, ,
\end{align*}
where $\P^{(1)}$ is the distribution of $\kappa(N)$ and $ \P^{(2)}$  the law of $\xi_{ij}$'s. For $\P^{(1)}$ a.e. $\omega_1$ 
$$
\lim_{N\to \infty}\P^{(2)}_{\kappa(N)(\omega_1)}\Big( \ \Big|\, \frac{ Y(t) }{ \rho^t U(\omega_1) N^{1- tr} } -1 \, \Big| \geq \varepsilon \ \Big) = 0\,,
$$
and we conclude by dominated convergence. \hfill $\Box$
\vspace{.2cm}

In \cite{Cook1989} Cook and Derrida consider the particular case where $\rho =1$ in (\ref{equa:defi.2.states}). From Corollary \ref{cor:fluct.Y.random}, we see that $Y(t)/N^{1-rt}$ converges in probability to one. Since under $\P_N$, $Y(t)$  is equal to the number of paths with zero energy at time $t$, the typical number of such paths is $N^{1-rt}$.





\subsection{Front Speed}\label{subsec.front.speed}

In this Subsection, we give the exact asymptotic for the front speed, proving Theorem \ref{teo:derrida}. The front positions can be computed by counting the number of times $Z$ visits $\oplus$. Indeed, at a given time $t$ either the front moves backwards and $\phi\big( X(t) \big) = \phi\big( X(t-1) \big) - 1$ or it stays still and $\phi\big( X(t) \big) = \phi \big( X(t-1) \big)$. We obtain that
$$
\frac{-N_t}{t} = \frac{\phi\big( X(t) \big) }{t},
$$
where $N_t$ is the stochastic process that counts the number of times that $Z$ visited $\oplus$ until time $t$. 

A classic result from renewal theory (see {\it e.g.} \cite{Durrett2010}) states that
\begin{equation}
\lim_{t \to \infty}\frac{N_t}{t} = \frac{1}{\E_\oplus[\tau]} \, .
\end{equation}
Hence, to determine the front velocity, it suffices to determine $\E_\oplus [  \tau ] \, . $
\begin{equation}\label{equa.integral.tau.sum.prob}
\E_\oplus [ \tau ] = \sum_{t = 0 }^{\infty} \P_\oplus ( \tau \geq t+1 ) =\sum_{t=0}^{\infty} \P_\oplus ( Y(t) \geq 1 ) \, .
\end{equation}
A consequence of Corollaries \ref{cor:Y.big.t}, \ref{cor:converg.Y.eta.random} and \ref{cor:fluct.Y.random} is that if $\xi$ is distributed according to (\ref{equa:defi.2.states}) with $\eta > 0$, then
$$
\lim_{N \to \infty} \P_\oplus \big(  Y(t) \geq 1  \big) = \left\{ 
\begin{array}{lcl}
1\, , &  \text{if}  & t \leq m  \, ; \\
0 \, , &   \text{if} & t \geq m+1  \,.
\end{array}
\right.
$$
Whereas we have the following limits when  $\eta = 0$
$$
\lim_{N \to \infty} \P_\oplus \big(  Y(t) \geq 1 \big) = \left\{ 
\begin{array}{lcl}
1 \, , &  \text{if}  & t \leq m-1  \, ; \\
1-e^{\rho^{m}} \, , &   \text{if} & t = m \, ; \\
0 \, , &   \text{if} &  t \geq m+1 \, .\\ 
\end{array}
\right.
$$
Then, to finish the proof of Theorem \ref{teo:derrida}, it suffices to show that
\begin{equation} 
\lim_{N\to\infty}\sum_{t\geq m+1} \P_\oplus ( Y( t) \geq 1 ) = 0\,. 
\end{equation}
Since $Y$ is a homogeneous Markov chain we use the Markov property at time $m+1$ to obtain
$$ 
\sum_{t \geq m +1 } \P_\oplus \big(  Y(t)\geq 1 \big) =\sum_{t=0}^{\infty} \sum_{k=1}^N \P_k \big( Y(t) \geq 1 \big) \,\P_\oplus \big(  Y(m+1) = k \big) \, .
$$
It is not difficult to see that under $\P_k$, $Y$  is stochastically dominated by $Y$ under $\P_N$, which implies that $\P_k \big(  Y(t) \geq 1 \big) \leq \P_N \big( Y(t) \geq 1  \big) \, $. Then,  applying this inequality in the above expression, we get 
\begin{equation}\label{equa:est.sum}
\sum_{t \geq m +1 } \P_\oplus \big(  Y(t) \geq 1  \big) \leq \P_\oplus \big( Y(m+1) \geq 1 \big) \E_\oplus[ \tau ] \, .
\end{equation}

\begin{prop}\label{prop:bound.tau} Let $\xi$ be distributed according to (\ref{equa:defi.2.states}). Then, $\E_x[\tau]$ is bounded in $N$
\begin{equation}
\sup_{N \in \N }\sup_{x\in\Omega(N)} \{\,\E_x[\tau ]\,\} < \infty
\end{equation} 
\end{prop}
\emph{Proof.} By Corollary \ref{cor:fluct.Y.random}, $\lim_{N \to \infty} \P_\oplus \big( \tau \geq m+2 \big) = 0 \,$. Therefore, there exists a constant $c_{\thec}<1$ \setcounter{aux}{\value{c}} such that for $N$ sufficiently large
$$
\P_\oplus ( \, \tau \geq m+2 \, ) \leq c_{\arabic{aux}} \, .
$$
Coupling the chains started from $\delta_x$ and $\delta_\oplus$ we obtain that $\P_{x}( \tau \geq m+2 ) \leq \P_{\oplus}( \tau \geq m+2 ) \, $ for every $x\in\Omega(N)$ and therefore
\begin{equation}\label{equa.prop:bound.tau}
\P_x ( \tau \geq m+2 ) \leq c_{\arabic{aux}} \, .
\end{equation}
Then, Proposition \ref{prop:bound.tau} follows as a consequence of the Markov property and (\ref{equa.prop:bound.tau}). In Subsection \ref{bound} we present an equivalent argument in all detail.\hfill $\Box$ 
\medskip

Applying Proposition \ref{prop:bound.tau} and Corollary \ref{cor:Y.big.t} in (\ref{equa:est.sum}), we conclude that
$$
\sum_{t \geq m +1 } \P_\oplus ( \,   Y(t) \geq 1 \,  ) = \mathcal{O} (\,N^{1 - (m+1) \, r }\,) \, .
$$
Hence, from (\ref{equa.integral.tau.sum.prob}) we obtain the limits  
\begin{equation} \label{equa:expected.value.tau.2.states}
\lim_{N \to \infty}\E_\oplus [\tau] = \left\{ 
\begin{array}{lcl}
1+ m , & \mbox{if} & r \neq 1/m \\ 
1+ m - e^{- \rho^{m}} , & \mbox{if} & r = 1/m \, ,
\end{array}
\right.
\end{equation}
proving  Theorem \ref{teo:derrida} in the case $r>0$.  

To finish the proof of Theorem \ref{teo:derrida} it remains to study the case where $r=0$. For that we use a coupling argument. Up to the end of this Subsection we denote by $ \xi (r)$
$$
\P( \xi (r) = 0 ) =1 - \P(\xi (r) = -1 ) \sim \rho/N^{1+r} \, .
$$
For $r>0$, the random variables $\xi(0)$ are stochastically larger than $\xi(r)$ for $N$ large enough. Denoting by $X_i^r (t)$ the stochastic process defined by $\xi(r)$ we construct the process in such a way that the following relation holds
$$
 0 \geq \frac{\phi\big( X^0(t) \big) }{t} \geq \frac{\phi\big( X^r(t) \big) }{t}.
$$
From (\ref{equa:expected.value.tau.2.states}), if we choose $r$ such that $1/r$ is not an integer, we have the lower bound
$$
0 \geq v_N (0) \geq v_N (r) \to \big( \  1+ \lfloor  1/ r \rfloor  \ \big)^{-1};
$$
whence taking $r$ to $0$, we have that $\lim v_N(0) = 0 $, which concludes the proof of Theorem  \ref{teo:derrida}.

\section{Front speed for the infinitely many states percolation distribution} \label{sec.front.speed.three.states}

In this Section, we prove a discrete version of Theorem \ref{teo.speed.3.states.gen}. We consider the case where $\xi_{ij}$ is defined as in (\ref{equa:defi.3.states}).

{\bf Assumption (A).} The random variable $\xi$ distributed according to (\ref{equa:defi.3.states}) satisfies Assumption (A) if there exists a constant $0 <\theta < 1 $ such that
$$
\lim_{N \to \infty} q_2 = \theta\,,
$$ 
and if $\vartheta $ defined in (\ref{equa.defi.vartheta}) is integrable.

\vspace{.2cm}   
In the non-critical case we do not need to assume the convergence of $q_2$. We prove Theorem \ref{teo.speed.3.states} under the weaker condition.

{\bf Assumption (A').} The random variable $\xi$ distributed according (\ref{equa:defi.3.states}) satisfies Assumption (A') if there exists a constant $0 <\theta' < 1$ such that for $N$ large enough
$$
\frac{q_2}{ ( 1-p_0 )} \leq \theta' \,,
$$ 
and if $\vartheta $ defined in (\ref{equa.defi.vartheta}) is integrable.
\vspace{.2cm}

We adapt the notation of the previous Section and let $Z(t) := ( Z_l (t)\,  ; l \in \N )$ be defined as
$$
Z_l(t) := \sharp \big\{ j \,;  1 \leq j \leq N \,, X_j(t) = \phi\big( X(t - 1) \big) - l \big\}\,.
$$
Then, $Z$ is a homogeneous Markov chain on the set 
$$
\Omega(N) := \bigg\{  x \in \{ 0, 1 ,  \ldots , N \}^{\N} \, ;  \sum_{i=0}^{\infty} x_i  = N  \bigg\},
$$ 
where $x_i$ are the coordinates of $x$. If at time $t$ we have that $Z(t)=x\in \Omega(N)$, it means that for each $k \in \N$ there are $x_k$ particles in position $-k\,$ with respect to the leader. In this situation, suppose that $x_0 \geq 1 $. Then the probability that at time $t+1$ there is some particle in position $-k$ with respect to the leader at time $t$ is given by,
\begin{equation}\label{equa:mult.coef}
s_k (x) := \bigg( \sum_{i=1}^{\infty} p_i \bigg)^{x_{k-1}} \ldots \bigg( \sum_{i=k}^{\infty} p_i \bigg)^{x_{0}} - \bigg( \sum_{i=1}^{\infty} p_i \bigg)^{x_{k}} \ldots \bigg( \sum_{i=k+1}^{\infty} p_i \bigg)^{x_{0}}\,,
\end{equation}
where we define $x_{-1}=0 $. So the probability that one particle has not yet moved at time $t+1$ is given by
$$ 
s_0(x) : = 1-\bigr( 1- p_0\bigl)^{x_0}\, .
$$
If $x_0 = 0$, we shift $( x_{0}, x_{1} ,\ldots )$ to get a nonzero first coordinate obtaining a vector $\tilde{x} \in \Omega(N)$ such that $\tilde{x}_0 \geq 1$. Then, one can check that 
$$
s_r (x) = s_r (\tilde{x})\,.
$$
The transition probability of the Markov chain $Z$ is given by
\begin{equation}\label{equa.Z.mult.nom.inf.range}
\P \big( Z(t+1) = y \mid Z(t)=x \big) = \mathcal{M}\big( N;s(x) \big) (y)\,,
\end{equation}
where $s(x) = \big( s_0(x), s_1(x) \ldots \big)$ and $\mathcal{M}\big( N;s(x) \big)$ denotes a Multinomial distribution with infinitely many classes. We refer to \cite{Comets2013}, Section 6, for more details on the computations. It is clear that $Z_0(t)$ has the same transition probability as the process studied in the two states model. In particular, the results proved in Subsection \ref{sec:number.particles.zero} hold with the obvious changes.

For a stopping time $T$, we define recursively $T^{(0)} = 0$ and for $i \geq 1$ 
\begin{equation}\label{equa.defi.T^i}
T^{(i)} (\omega) : = \inf \{ t > T^{(i-1)} ( \omega ) ;   t = T \circ \Theta_{T^{(i-1)} ( \omega )} ( \omega ) \} \,,
\end{equation}
where $\Theta_t$ is the time-shift operator. We adopt the convention that $\inf \{ \, \emptyset \,\} = \infty$. Once more we denote by $\tau$ the stopping time defined as
\begin{equation}\label{equa.defi.tau.three.states}
\tau := \inf \big\{ t\,;\,\phi\big(X(t)\big) < \phi\big(X(t-1)\big) \big\} \,.
\end{equation}
In contrast with the previous Subsection, $\tau$ is not a renewal time for $Z$. Let $T_x$ be the first time that $ Z(t)$ visits $x$
\begin{equation}
T_{x}: = \inf \{ t \,  ; Z(t) = x \} \,.
\end{equation}
We adapt the notation of Section \ref{sec.derrida.two} and define $\oplus := \big(N,0, \ldots \big) \in \Omega(N)$ and $\triangle : =  (0,N,0,\ldots)  \in \Omega(N)\,.$ Finally, we keep notation (\ref{eq:integer.part.r}) and let $m$ be the integer part of $1/r$ and $\eta$ its fractional part. 

We now state the main result of the Section.

\begin{teo} \label{teo.speed.3.states} Assume that $\xi$ satisfies Assumption (A). Then 
\begin{equation}
\lim_{N \to \infty} v_N = \left\{ 
\begin{array}{lcl}
-\big( \,1+ \lfloor  1/r  \rfloor  \big)^{-1}, & \text{if} & 1/r \not\in \N   \\
- \big( \, \lfloor   1/r  \rfloor  + 1- 1/g(\theta) \, \big)^{-1} , & \text{if} & 1/r = m \in \N \,,
\end{array}
\right.
\end{equation}
where $g (\theta ) \geq 1 $ is a non-increasing function. The conclusion in the case $r \neq 1/m$ still holds if $\xi$ satisfies the weaker Assumption (A').
\end{teo}

\subsection{The Distribution of $Z(\tau)$}\label{sec:jumping}

In this Section we study the limit distribution of $Z(\tau)$ as $N \to \infty$. When $\eta>0$ the limit is similar to the one obtained in the previous results. 

\begin{prop} \label{prop:conv.prob.Xtau} Assume that $\xi$ satisfies Assumption (A') and that $\eta>0$. Then,
\begin{equation}
\lim_{N \to \infty } \P_\oplus \big(  Z(\tau) = \triangle   \big) = 1   \, .
\end{equation}
\end{prop}
\vspace{0.2cm}

The case $\eta = 0$ is critical. We show that $Z_1(\tau)/N$ converges in distribution and that the limit distribution is a functional of a Poisson random variable. 

\begin{prop} \label{prop:conv.z.tau.under.oplus} Assume that $\xi$ satisfies Assumption (A) with $\eta=0$ . Then under $\P_\oplus$, $Z_0 (m)$ converges in distribution to $\Pi(\rho^m)$ a Poisson random variable with parameter $\rho^m$.

Moreover, there exists a function $\,G\,:\, \N \to [0 , 1]$ (see Definition \ref{equa:G}) such that
\begin{equation}
\bigg( \frac{Z_1(\tau)}{N}\,, \sum_{i=2}^{\infty} \frac{Z_i(\tau)}{N} \bigg) \stackrel{d}{\rightarrow} \Big(G\big(\Pi(\rho^m)\big), 1-G\big(\Pi(\rho^m)\big) \Big)\,.
\end{equation}
\end{prop}

Before analyzing the cases $\eta = 0 $ and $\eta >0$ separately, we prove a technical Lemma that holds in both cases. It can be interpreted as follows: if at time $t$ there are sufficiently many leading particles, then at $t+1$, with high probability, there is no particle at distance two or more from the leaders at $t$.

\begin{lem} \label{lem:techinical.multi} Assume that $\xi$ satisfies Assumption (A'). For $x = x(N) \in \Omega(N)$ such that  
$$
\log N = o ( x_0 ) \, ,
$$
define $ s_{i}(x) $ as in (\ref{equa:mult.coef}) and let $\mathcal{M}\big( N;s(x) \big)$ be a Multinomial random variable with infinitely many classes as in (\ref{equa.Z.mult.nom.inf.range}). Then,
\begin{equation}
\lim_{N\to \infty}\P \bigg( \mathcal{M}\big( N;s(x) \big) \in  \Big\{ y \in \Omega(N) \, ;  \sum_{i=2}^{\infty} y_{i} = 0  \Big\} \, \bigg) = 1 \,  .
\end{equation}
\end{lem}

\emph{Proof.} We can write
\begin{align*}
\P \bigg( & \mathcal{M}\big( N;s(x) \big) \in  \Big\{ y \in \Omega(N) \, ;  \sum_{i=2}^{\infty} y_{i} = 0  \Big\} \, \bigg) \\
 & = \sum_{n=0}^{N} \P \bigg( \mathcal{M}\big( N;s(x) \big) \in  \Big\{ y \in \Omega(N) \, ; y_{0} = n \,, y_1 =N-n  \Big\} \, \bigg) \\
&=  \sum_{n=0}^{N} \frac{N!}{n!(N-n)!}s_0(x)^n \, s_{1}(x)^{N-n} 
\geq \big( 1 - \theta'^{\,x_0} \big)^N \, ,
\end{align*}
where the last inequality holds for $N$ large enough as a consequence of Assumption (A'). Since $o( x_0 ) = \log N$ we obtain that $\big( 1 - \theta'^{\,x_0} \big)^N \to 1 $, proving the result. \hfill $\Box$
\vspace{0.2cm}


\subsubsection*{Case $\eta>0$}

We have already introduced all necessary tools to prove Proposition \ref{prop:conv.prob.Xtau}.

\emph{Proof of Proposition \ref{prop:conv.prob.Xtau}. } From Corollaries \ref{cor:Y.big.t} and \ref{cor:fluct.Y.random} we see that $\P_\oplus \big( \tau \neq m+1 \big) \to 0 \,.$ Then, it suffices to prove that $\P_\oplus \big(  Z(m+1) = \triangle   ;  \tau = m+1 \big) \to 1 \,.$ 
$$
\P_\oplus  \big( Z(\tau ) = \triangle \, ;  \tau = m+1  \big)= \sum_{x \in \Omega(N)} \P_\oplus \big(  Z(m+1)  = \triangle \, ;Z(m) = x \, ; \tau = m+1 \big) \, .
$$
Since $\tau = m+1$ it suffices to consider $x$ such that $x_0 \geq 1$. Fix $0<\varepsilon< \rho^{ m}\,$ and take $x\in \Omega(N)$ such that $| x_0/N^{r \, \eta} -\rho^{m} | < \varepsilon$. From (\ref{equa.Z.mult.nom.inf.range}),  
\begin{align}\label{equa.proof.teo.Z.to.triangle}
 \P_\oplus & \big(  Z_0(m+1) = \triangle | Z_0 (m) = x \big) \nonumber \\
&= \mathcal{M}\big( N;s(x) \big) ( \triangle )  = s_{1}(x)^N \nonumber\\
 &= \Big( \big( 1- p_0\big)^{x_0}  - \big(1- p_0\big)^{x_{1}} \big(1- p_0 - p_1\big)^{x_0}\Big)^N \nonumber \\
 & \geq \big(1- p_0\big)^{x_0 N} \big(1  -  \theta'^{\,x_0}\big)^N \, ,
\end{align}
where the last inequality is a consequence of Assumption (A'). Since $x_0 = \mathcal{O} (\, N^{\eta r} \,)$ we conclude after a first order expansion that the lower bound in (\ref{equa.proof.teo.Z.to.triangle}) converges to one. Moreover, the rate of convergence is bounded from below by
$$
\big(1- p_0\big)^{(\rho^m+\varepsilon ) N^{1+r \eta}} \big(1  -  \theta'^{(\rho^m-\varepsilon ) N^{r \eta}}\big)^N \,,
$$
which converges to one. Then, by Proposition \ref{prop:large.deviation.Y} and Equation (\ref{equa.proof.teo.Z.to.triangle}), we see that
$$
\P_\oplus  \big(  Z(\tau) = \triangle \big) \geq \sum_{| x_0 / N^{r  \eta} -\rho^{m} | < \varepsilon } \P_\oplus  \big(  Z(\tau ) = \triangle ;  Z(m) = x ; \tau = m+1 \big) \\
$$
converges to one, proving the result. \hfill $\Box$ 

\subsubsection*{Case $\eta=0$} \label{sec:jumping.eta.zero} 

In this paragraph, we prove Proposition \ref{prop:conv.z.tau.under.oplus} and also a generalization that allows us to compute the distribution of $Z_1(\tau^{(i)})$.



\begin{lem} \label{lem:mult.techinique.many.particles}  Assume that $\xi$ satisfies Assumption (A') with $\eta=0$. Fix $ \, 0<a<b \, $ and denote by $\Omega_a^b(N)$ the subset of $\Omega(N)$ defined as
$$
\Omega_a^b(N):= \big\{ x \in \Omega(N) \, ;  a N^{1/m} \leq x_0  \leq b N^{1/m}   \big\}\,.
$$
Then the following limit holds
\begin{equation}
\lim_{N \to \infty} \,\,\, \sup_{x \in \Omega_a^b (N)} \P_x \big( Z(1) \neq \triangle \mid Z_0(1) = 0 \big) = 0 \,.
\end{equation}
\end{lem} 
\emph{Proof.}  It is not difficult to obtain the following inequality
$$
\P_x \big( Z(1) \neq \triangle \mid Z_0(1) = 0 \big) \leq  \frac{ \P_x \Big( Z(1) \in \big\{ y \in \Omega(N) \, ;  \sum_{i=2}^{\infty} y_{i} \neq 0 \big\} \Big) }{\P_x \big( Z_0(1) = 0  \big)} \, .
$$ 
From (\ref{equa.Z.mult.nom.inf.range}) we have that under $\P_x$, $Z(1)$  is distributed according to $\mathcal{M}\big(N,s(x)\big)$ . Then, as a consequence of Lemma \ref{lem:techinical.multi} 
\begin{align*}
\P_x & \bigg( Z(1) \in \Big\{ y \in \Omega(N) \, ;  \sum_{i=2}^{\infty} y_{i} \neq 0 \Big\} \bigg)\\
 &= 1-\P \bigg( \mathcal{M}\big( N;s(x) \big) \in  \Big\{ y \in \Omega(N) \, ;  \sum_{i=2}^{\infty} y_{i} = 0  \Big\} \, \bigg) \to 0 \,.
\end{align*}
Moreover, the rate of decay is bounded from above by 
$$
1- \big(1-\theta'^{a N^{1/m}}\big)^N \to 0 \,.
$$

To finish the proof it suffices to show that $ \P_x \big( Z_0(1) = 0  \big) $ is bounded away from zero. Indeed,under $\P_x$, $Z_0 (1)$  is distributed according to a Binomial random variable of parameter $N$ and $s_0(x)$. 

From the hypotheses of the Lemma we have that 
$$
s_0(x) \geq 1-(1-p_0)^{b\,N^{1/m}}\,.
$$
Coupling $Z(1)$ with $\mathcal{B}$ a Binomial of parameter $N$ and $1-(1-p_0)^{b\,N^{1/m}}$, we conclude that
$$
 \P_x \big( Z_0(1) = 0  \big) \geq \mathcal{B} \big( N, 1-(1-p_0)^{b\,N^{1/m}}\big)(0)   \to e^{-\rho b} \,.
$$
\hfill $\Box$ 
\vspace{.2cm}

From Corollary \ref{cor:fluct.Y.random}, we see that under $\P_\oplus$, $Z_0 (m-1) /N^{1/m}$  converges in probability to $\rho^{m-1}$, as $N\to\infty$. Hence, from Lemma \ref{lem:mult.techinique.many.particles}, we conclude that  
\begin{equation}\label{equation.convergence.Z(m)=0}
\lim_{N \to \infty } \P_\oplus \big( Z(\tau) = \triangle  \,| \, Z_0(\,m\,) = 0 \big) = 1 \,.
\end{equation}
This is the first step to prove Proposition \ref{prop:conv.z.tau.under.oplus}. The second step is to study the conditional distribution of $Z(\tau)$ under $Z_0(m) = x_0$ for a positive integer $x_0$. 




\begin{prop}\label{prop:conv.of.Z(1)} Assume that $\xi$ satisfies Assumption (A) with $\eta=0$. Let $ k $ be a nonzero integer and denote by $\Omega_k(N)$ the subset of $\Omega(N)$ defined as
$$
\Omega_k(N) : = \big\{ x \in \Omega(N) \, ; x_0 = k  \big\} \,.
$$
Then, for $\varepsilon >0 $ the following limit holds 
\begin{equation}
\lim_{N \to \infty} \, \sup_{ x \in \Omega_k(N) } \, \P_x \bigg( \,\, \bigg| \bigg ( \frac{Z_1 (1)}{N} \, ,  \frac{ \sum_{i\geq 2}Z_i (1)}{N} \bigg) - (1-\theta^k,  \theta^k)   \bigg| > \varepsilon  \bigg) = 0 \, .
\end{equation}

\end{prop}

\emph{Proof.} From (\ref{equa.Z.mult.nom.inf.range}), we see that under $\P_x$, $Z(1)$ is distributed according to an infinite range Multinomial of parameters $N$, and $s(x)$. In particular under $\P_x$ the triplet
\begin{center}
$
\big( Z_0 (1)\,, Z_1 (1) \, , \sum_{i\geq 2}Z_i (1) \big) \,,
$
\end{center}
is distributed according to a three classes Multinomial of parameters $N$ and $ \big( s_{0}(x) , s_{1}(x) , \sum s_i(x) \big)$. If $\xi$ satisfies Assumption (A) and $x\in\Omega_k(N)$, we have that 
\begin{equation}
\lim_{N\to \infty} s_0(x)  =0 \,;\quad \lim_{N\to \infty} s_{1}(x)  = 1-\theta^k \,; \quad \lim_{N\to \infty} \sum s_{i}(x)  = \theta^k \, .
\end{equation}
The rate of convergence is uniform on $x \in \Omega_k (N)$. 

A three classes Multinomial random variable as above satisfies a large deviation principle (see {\it e.g.} \cite{Dembo2010, DenHollender2008}) and the rate function is given by
\begin{equation}\label{equa:mult:larg.dev}
\Lambda^*(y) =  \left\{
\begin{array}{ccl}
y_{1} \log \Big( \frac{(\theta^k) y_{1} }{(1 -y_{1})( 1-\theta^k )}  \Big) - \log \Big( \frac{\theta^k  }{1 -y_{1}} \Big) & \mbox{ if } &  y_{1} + y_{2} = 1 \, , \\
\infty  &  &  \text{  otherwise.  }
\end{array}
\right.
\end{equation}
The only zero of $\Lambda^*$ is at $y = (  0,  1-\theta^k, \theta^k )$. Implying the convergence in probability
\begin{center}
$
\frac{1}{N}\big( Z_0 (1)\,, Z_1 (1) \, , \sum_{i\geq 2}Z_i (1) \big) \to (  0,  1-\theta^k, \theta^k ) \,.
$
\end{center}   \hfill $\Box$
\vspace{.2cm}

We now give the definition of the function $G$ appearing in Proposition \ref{prop:conv.z.tau.under.oplus}. 

\begin{defi}
Let $G: \N \longrightarrow [0,1] $ be defined as
\begin{equation}\label{equa:G}
G(k) =  \left\{
\begin{array}{ccl}
1 - \theta^k, & \mbox{ if } & \, k \geq 1 \\
1 , & \mbox{ if } & \, k =0  \\
\end{array}
\right.
\end{equation} 
where $\theta$ is given by Assumption (A).
\end{defi} 

\emph{Proof of Proposition \ref{prop:conv.z.tau.under.oplus}.\,} From Corollary \ref{cor:converg.Y.eta.random}, we have that under $\P_\oplus$, $Z_0 (m)$ converges in distribution to a Poisson random variable of parameter $\rho^m$. Hence, to prove Proposition \ref{prop:conv.z.tau.under.oplus} it suffices to show that
\begin{align*} \label{equa:teo.favorite.site.to.jump.eta.0}
\P_\oplus  & \bigg( \,\, \bigg| \bigg ( \frac{Z_1 (\tau)}{N} \, ,  \frac{ \sum_{i\geq 2}Z_i (\tau)}{N} \bigg) - \Big( G \big( Z_0(m) \big),  1- G \big( Z_0(m) \big) \Big)   \bigg| > \varepsilon  \bigg) \\
 &= \sum_{k= 0}^{N}  \P_\oplus \bigg( \,\, \bigg| \bigg ( \frac{Z_1 (\tau)}{N} \, ,  \frac{ \sum_{i\geq 2}Z_i (\tau)}{N} \bigg) -  \big( G(k),  1- G (k) \big)   \bigg| > \varepsilon \,; Z_0(m) = k  \bigg)\,
\end{align*}
converges to zero. From (\ref{equation.convergence.Z(m)=0}) and Proposition \ref{prop:conv.of.Z(1)}, we know that for each $k \in \N$ the terms in the above sum converge to zero. Then, from the tightness of $Z_0 (m)$ we obtain that the sum itself converges to zero, proving the result.\hfill $\Box$

\medskip

We finish the present Subsection by computing the limit distribution of $Z(\tau^{(i)})$ for $i \in \N$. We also prove the convergence of some related processes that will appear when calculating the front velocity in Subsection \ref{sec:front.speed.3.states}. 
\begin{prop}\label{prop:chain.started.random.position} Assume that $\xi$ satisfies Assumption (A) and that $\eta=0$. Let $\kappa(N)$ be random variables taking values in $\Omega(N)$, such that $\kappa(N)$ and the $\xi_{ij}$ are independent for every $N$. Denoting by $\kappa_0(N)$ the first coordinate of $\kappa(N)$ we also assume that
\begin{equation}\label{nontrivial.cond}
\lim_{N \to \infty} \frac{\kappa_0(N)}{N} = U \, \, \, a.s.
\end{equation}
where $U$ is a positive random variable.  Then, under $\P_{\kappa(N)}$, we have that 
\begin{enumerate}
\item $Z_0( m )$ converges weakly to $V$ a doubly stochastic Poisson random variable which distribution is determined by the Laplace transform
\begin{equation}
\E[  e^{s V}  ] = \E \big[ \exp(  e^s-1  )\,\rho^{m} \, U \big] \, .
\end{equation}
\item Furthermore, the joint convergence also holds 
\begin{equation} \label{equa:conv.Z.lambda-1.Z.tau}
\bigg( Z_0( m), \frac{Z_{1} ( \tau )}{N} , \tau \bigg) \stackrel{d}{\longrightarrow } \big( V  ,G(V), m+1_{\{V \neq 0\}} \big) \,.
\end{equation} 
\end{enumerate}
\end{prop}

\emph{Proof.} Since $ \kappa_0(N)/ N \to U $ the hypotheses of Corollary \ref{cor:converg.Y.eta.random} are satisfied, implying the first statement of the Proposition. From Corollaries \ref{cor:Y.big.t} and \ref{cor:fluct.Y.random}, we see that $\P( m \leq \tau \leq m+1)$ converges to one. Since $\tau = m $ if and only if $Z_0(m) = 0$ we obtain the convergence in distribution
$$
\tau \stackrel{d}{\to} m+1_{\{V \neq 0\}} \,.
$$  
Finally, to prove that $Z_{1} (\tau)/N$ converges to $G(V)$, we proceed as in Proposition \ref{prop:conv.z.tau.under.oplus} and show by dominated convergence that
$$
\lim_{N\to \infty} \E\Bigg[ \P_{\kappa(N)} \bigg( \,\, \bigg| \bigg ( \frac{Z_1 (\tau)}{N} \, ,  \frac{ \sum_{i\geq 2}Z_i (\tau)}{N} \bigg) - \big( G \big( Z_0(m) \big),  1- G \big( Z_0(m) \big) \big)   \bigg| > \varepsilon  \bigg) \Bigg] =0 \,.
$$
 \hfill $\Box$
\vspace{.2cm}

As an application of Proposition \ref{prop:chain.started.random.position} we can calculate the distribution of $Z(\tau^{(2)})$. Indeed we can consider the convergence in Proposition \ref{prop:conv.z.tau.under.oplus} as the stronger a.s. convergence. We do not lose any generality since we can construct a sequence of random variables (possibly in an enlarged probability space) $ \kappa(N) \stackrel{d}{=}Z(\tau)$ that converges a.s. Details about this construction can be found in \cite{Billingsley2009}. Passing to the appropriate product space we also consider that the $\kappa$'s and $\xi_{ij}$'s are independent. Then by the strong Markov property, we obtain that
\begin{equation}
\P_{\kappa (N)} \big( Z \circ \Theta_{\tau} (t) \in \cdot \big) \stackrel{d}{=}\P_{Z(\tau)} \big(Z \circ \Theta_{\tau} (t) \in \cdot \big) \,,
\end{equation}
for $t \geq 0 \,$. Then, under $\P_\oplus$ we obtain that
$$
\Big(Z_0( \tau + m  ) \, , \, \frac{Z_{1} (\tau^{(2)})}{N}  \,, \tau^{(2)} - \tau^{(1)} \Big)  \stackrel{d}{\longrightarrow } \big ( V^{(2)} , G(V^{(2)}), m+1_{\{V^{(2)} \neq 0\}}\big) \,,
$$
where $V^{(2)}$ is a doubly stochastic Poisson variable governed by $V^{(1)}$ the limit distribution of $Z(m)$. This method can be iterated to obtain the following result.

\begin{lem}\label{lem:joint.conv.finite.dem.proj} Assume that $\xi$ satisfies Assumption (A) with $\eta=0$. Let $\Delta \tau^{(i)}_N = \tau^{(i)} - \tau^{(i-1)} $. Then, under $\P_\oplus$, 
\begin{equation}
\big\{ \big(  Z_0 (\tau^{(i-1)} + m) \,, Z_{1}(\tau^{(i)})/N  \, , \Delta \tau^{(i)}_N   \big) \,; 1\leq i \leq l \big\} \,, 
\end{equation} 
converges weakly to
\begin{equation}
\Big\{ \, \Big( V^{(i)} , G( V^{(i)} ) \,,  m+ 1_{\{ V^{(i)} \neq 0 \}} \Big) \,; 1\leq i \leq l \Big\} \,.
\end{equation}
where $l \in \N$ is fixed. The distribution of $V^{(i)}$ are determined by
\begin{align}
\P & \big( V^{(i+1)} = l \, \mid \, V^{(j)} =  t_j , j \leq i \big)  \nonumber \\
& = \P \big(  V^{(i+1)} = l \, \mid \, V^{(i)} =t_i  \big) = e^{-G(t_i) \rho^{m}}\frac{( G(t_i) \rho^{m}  )^l}{l!} \, ,
\end{align} 
where $2 \leq i \leq k-1\, ,t_j \in \N$ and $V^{(1)}$ is distributed according to a Poisson variable with parameter $\rho^{m}$.
\end{lem}

\emph{Proof.} It is a direct consequence of Proposition \ref{prop:chain.started.random.position} and an induction argument. \hfill $\Box$
\vspace{.2cm}

With a very small effort we can state Lemma \ref{lem:joint.conv.finite.dem.proj} in a more general framework. We consider the space of real valued sequences $\R^\N$ where we define the metric
$$
d(a,b) = \sum_{n=0}^{\infty} \frac{| a_n - b_n |}{2^n} \,.
$$
A complete description of this topological space can be found in \cite{Billingsley2009}. Since time is discrete, the following Proposition holds as a Corollary of Lemma \ref{lem:joint.conv.finite.dem.proj}.

\begin{prop}\label{prop:conv.joint.law.process} Assume that $\xi$ satisfies Assumption (A) with $\eta=0\,$. Then, under $\P_\oplus$ the process
\begin{equation}
\big\{ \, \big( Z_0(\tau^{(i-1)} +m)  \, , Z_{1}( \tau^{(i)} )/N \, , \Delta \tau^{(i)} \big) ; i \in \N \big\}
\end{equation}
converges weakly in $\big(  (\,\R^\N \,)^3  ,d \big)$. The limit distribution $\mathbb{W}_\theta$ is given by
\begin{equation}
\big\{ \, \big( V^{(i)},  G(V^{(i)}) ,   m+1_{\{ V^{(i)} \neq 0 \}}  \big) ; i \in \N  \big\}\,,
\end{equation}
where $V^{(i)}$ is a Markov chain with initial position at $0$ and transition matrix given by
\begin{equation}
\P \big( V^{(i+1)} = l \mid  V^{(j)} = t \big) = e^{-G(t) \rho^{m}}\frac{(G(t) \rho^{m}  )^l}{l!} \, ,
\end{equation} 
that is a Poisson distribution with parameter $\rho^m G(t)$.  
\end{prop} 

\subsubsection*{Process convergence in the case $\eta>0$}

For the sake of completeness, we state the result in the case $\eta>0$. We omit the proof of the Proposition and leave the details to the reader . 

\begin{prop}\label{prop:conv.joint.law.process.eta.0} Assume that $\xi$ satisfies Assumption (A') and that $\eta>0\,$. Then under $\P_\oplus$ the process
\begin{equation}
\big\{  \big(  Z_{1}(  \tau^{(i)})/N  , \Delta \tau^{(i)} \big) ;  i \in \N  \big\}
\end{equation}
converges weakly in $\big( (\R^\N )^2  , d \big)$. The limit distribution is non-random, and concentrated on the sequence
\begin{equation}
\big\{ \, \big(  a_i  ,   b_i \big);  \quad a_i = 1 \mbox{ and } b_i = m+1 \  \forall i \in \N  \big\}\,.
\end{equation}
\end{prop}

\subsection{Uniform integrability and bounds for $T_\triangle$ } \label{bound}

In this Subsection, we show that if $\xi$ satisfies Assumption (A'), then $\E_x [ T_\triangle ]$ is bounded independently from the initial configuration $x$.
\begin{equation}\label{equa.aim.section.bound}
\sup_{N \in \N} \ \sup_{x \in \Omega(N)} \E_{x} [T_\triangle] < \infty \,.
\end{equation}
We prove (\ref{equa.aim.section.bound}) through the following steps.  
\begin{enumerate}
\item There exists a set, which we denote by $\Xi $, such that for $N$ large enough and every starting point $x \in \Xi$ there is a positive probability to visit $\triangle$ before $m+1$ 
\begin{equation} \label{equa:first.step.x.in.Lambda}
\P_x \big( T_\triangle \leq m+1 \big) > c_{\thec}  \,,
\end{equation}
\setcounter{aux1}{\value{c}} where $c_{\arabic{aux1}}>0$ does not depend on $x \in \Xi$ 
\item For $N$ sufficiently large and every starting point $x \in \Omega(N)$ there is a positive probability to visit $\Xi$ before $m+1$
\begin{equation} \label{equa:second.step.T.Lambda}
\P_x \big(  T_\Xi \leq m+1  \big)>c_{\thec} \, ,
\end{equation}
\setcounter{aux2}{\value{c}}where $T_{\Xi}:=\inf\{ t ; Z(t) \in \Xi \}$ and $c_{\arabic{aux2}} $ does not depend on $x \in \Omega(N)\,.$ 
\end{enumerate}
Before proving this two statements, we show that they indeed imply (\ref{equa.aim.section.bound}).

\begin{prop} \label{prop:bound.integral.T} Assume that $\xi$ satisfies Assumption (A'). Then
\begin{equation}\label{equa.prop:bound.integral.T}
\sup_{N \in \N}\sup_{x \in \Omega(N)} \E_{x} [T_\triangle] < K \,,
\end{equation}
where $K < \infty \,.$
\end{prop}
\emph{Proof.} If (\ref{equa:first.step.x.in.Lambda}) and (\ref{equa:second.step.T.Lambda}) hold, then for $N$ large enough and any starting point $x \in \Omega(N)$ 

\begin{align*}
\P_x & \big( T_\triangle \leq 2 m+2 \big) \\
& \geq \P_x \big( T_\triangle \leq 2 m+2; T_\Xi \leq m+1  \big) \\
& \geq \P_x \big( T_\triangle - T_\Xi \leq m+1 ; T_\Xi \leq m+1  \big) \\
&= \E_x \big[ \P_{Z(T_\Xi )} [T_\triangle \leq m+1] 1_{ T_\Xi \leq m+1}  \big] \hspace{1.0cm}  \mbox{ (Markov property)} \\
& \geq c_{\arabic{aux1}} \, c_{\arabic{aux2}} >0\,.
\end{align*}
Let $c_{\thec} = 1 - c_{\arabic{aux1}} \, c_{\arabic{aux2}} < 1\,.$ \setcounter{aux}{\value{c}} Then it is clear that $\sup_{y \in \Omega (N)} \P_y (T_\triangle \geq 2 m + 3 ) \leq c_{\arabic{aux}} \,$. For $i \in \N$, let $j$ be such that $(2 m+3)j \leq i <(2 m+3)(j+1)$. Using the Markov property $j$ times we obtain the upper bound 
$$
\P_x (T_\triangle \geq i) \leq \biggl( \,\, \sup_{y \in \Omega(N)}  \bigl\{ \P_y \big( T_\triangle \geq (2 m+3) \big) \bigr\} \biggr)^j\,. 
$$
We now show that the expected value of $T_\triangle$ is bounded.
\begin{align*}
\E_x [T_\triangle] & = \sum_{i=0}^{\infty} \P_x(\,T_\triangle \geq i\,) \\
& \leq \sum_{j=0}^{\infty}(2 m+3)  \sup_{y \in \Omega(N)} \bigl\{ \P_y \big( T_\triangle \geq (2 m+3) \big) \bigr\}^j \\
& \leq \sum_{j=0}^{\infty} (2 m+3) c_{\arabic{aux}}^{\,j}  = \frac{(2 m+3)}{1-c_{\arabic{aux}}}.
\end{align*}
Therefore (\ref{equa.prop:bound.integral.T}) holds with $ K= (2 m+3) / ( 1-c_{\arabic{aux}} ) \,.$ \hfill $\Box$
\vspace{.2cm}

We now present the formal definition of $\Xi$.
\begin{defi} \label{defi:Lambda} For $x\in \Omega(N)$ define $I(x) = \inf\{ i \in \N ; x_i \geq 1 \}$. Then, $\Xi$ is the subset of $\Omega(N)$ defined as follows. 
$$
\Xi := \bigl\{ x \in \Omega(N)\, ; x_{I(x)} \geq \alpha N  \big\},
$$
where $0<\alpha<1- \theta'$ and $\theta'$ is given by Assumption (A'). Hence, if $Z(t) \in \Xi$ there are at least $ \alpha N$ leaders at time $t\,.$
\end{defi}

We prove (\ref{equa:first.step.x.in.Lambda}) and (\ref{equa:second.step.T.Lambda}) in the next two Lemmas.

\begin{lem} \label{lem:jump.good.set} Assume that $\xi$ satisfies Assumption (A'). Then, for $\Xi$ given by Definition \ref{defi:Lambda} there exists a positive constant $c_{\arabic{aux1}}$ such that for $N$ sufficiently large (\ref{equa:first.step.x.in.Lambda}) holds $\forall \, x \in \Xi$  .
\end{lem} 
\emph{Proof.} Note that
\begin{align*}
\P_x \big( T_\triangle \leq m+1 \big) & \geq \P_x \big( Z(\tau) = \triangle \, ; \tau \leq m+1 \big) \\
& = \P_x (  Z(\tau) = \triangle  ) - \P_x \big( Z(\tau) = \triangle \, ; \tau \geq m+2 \big)  \, ,
\end{align*}
From Corollary \ref{cor:Y.big.t}, the second term in the lower bound converges to zero as $N \to \infty$ and the rate of decay is uniform on $y \in \Omega(N)$. Hence it suffices to show that there exists a positive constant $c_{\thec}$  \setcounter{aux}{\value{c}} such that uniformly on $x \in \Xi $ 
\begin{equation}\label{lem:eq2}
\lim_{N \to \infty} \P_x ( \, Z(\tau) = \triangle \, ) \geq c_{\arabic{aux}} \,, 
\end{equation}
To prove (\ref{lem:eq2}) we distinguish between the cases $\eta = 0$ and $\eta >0$. We start with the latter case $\eta>0\,.$ Let $Y(t)= Z_0 (t)1_{\{t \leq \tau\}} $ and denote by $Y_k$ the process started from $\delta_k$. Then, for $x \in \Xi$ we can couple the processes in such a way that 
\begin{equation}\label{equa:sto.bounded}
 Y_{\lfloor \alpha N \rfloor} (t) \leq Y_{x_{I(x)}} (t) \leq Y_N (t),
\end{equation}
where $x_{I(x)}$ is the number of leaders when $Z(0)= x$. From the proof of Corollary \ref{cor:fluct.Y.random} and (\ref{equa:sto.bounded}) we obtain
\begin{equation}\label{equa:lem:many:Y}
\lim_{N \to \infty}\P_{x} \big( ( \rho^{m} - \varepsilon ) \alpha N^{\eta r} \leq Z_0 (m) \leq ( \rho^{m} + \varepsilon ) N^{\eta r} \big) = 1 \, .
\end{equation}
Finally, applying the arguments of Lemma \ref{lem:techinical.multi},
$$
\lim_{N \to \infty} \P_x ( \, Z (m+1 ) = \triangle \, ) =1 \,.
$$
In particular, any $0 < c_{\arabic{aux}} < 1 $ satisfies (\ref{lem:eq2}) for $N$ sufficiently large. 

The case where $\eta=0$ is similar but it requires an additional step. (\ref{equa:sto.bounded}) still holds, hence by the same arguments we obtain
$$
\P_{x} \Big( \, ( \, \rho^{m-1} - \varepsilon \, ) a N^{1/m} \leq Z_0 (m-1) \leq ( \, \rho^{m-1} + \varepsilon \, ) N^{1/m} \ \Big) = 1 \, .
$$
From Lemma \ref{lem:mult.techinique.many.particles}, we see that $\P_x \big( Z(\tau) = \triangle \mid Z_0 (m) = 0  \big) \to 1\,, $ and from the coupling argument (\ref{equa:sto.bounded}) and Corollary \ref{cor:converg.Y.eta.random} we obtain the following limit
$$
\P_{x} \big(  Z_0 (m) = 0 \big) \geq P_{\oplus} \big( Z_0 (m) = 0 \big) \to 1-e^{\rho^m} \, .
$$ 
Then, any $c_{\arabic{aux}} $ smaller than $ 1-e^{\rho^m}  $ satisfies (\ref{lem:eq2}) for $N$ sufficiently large, proving the statement. \hfill $\Box$ 
\medskip

\begin{lem} \label{lem:jumping.Lambda} Assume that $\xi$ satisfies Assumption (A'). Then, for $\Xi$ given by Definition \ref{defi:Lambda} there exists a positive constant $c_{\arabic{aux2}}$ such that for $N$ large enough (\ref{equa:second.step.T.Lambda}) holds $\forall \, x \in \Omega(N)$ .
\end{lem}

\emph{Proof.} Since $ \P_x ( \tau \geq m+2 )$ converges to zero uniformly on $x \in \Omega(N)$, it is sufficient to show that for $N$ sufficiently large  
$$
\P_x \big( Z ( \tau ) \in \Xi  \big) \geq c_{\thec}\,,
$$
\setcounter{aux}{\value{c}}where $c_{\arabic{aux}}>0$ does not depend on $x\in\Omega(N)\,$.   
\begin{align}\label{inq:4}
\P_x &( Z ( \tau ) \in \Xi  ) = \sum_{k = 1}^{\infty} \P_x \big( Z ( k ) \in \Xi ; \tau = k   \big) \nonumber \\
& = \sum_{k = 1}^{\infty} \sum_{y \in \Omega(N)}  \E_x \big[ \P_y ( Z ( 1 ) \in \Xi ; \tau = 1 ) 1_{ \{ Z(k-1)=y ;\, \tau \geq k \} }   \big] \quad \mbox{(Markov property)} \nonumber \\
& \geq \inf_{y \in \Omega(N)} \big\{  \P_y ( Z ( 1 ) \in \Xi \mid \tau = 1 ) \big\}  \sum_{k = 1}^{\infty} \sum_{y \in \Omega(N)}  \E_x \Big[  \P_y(\tau=1)  1_{ \{ Z(k-1)=y ; \, \tau \geq k \} }   \Big] \nonumber \\
& = \inf_{y \in \Omega(N)} \big\{  \P_y ( Z ( 1 ) \in \Xi \mid \tau = 1 ) \big\} \,.
\end{align}
Then, it suffices to show that the infimum in (\ref{inq:4}) is larger than $ c_{\arabic{aux}}$. Recall that under $\P_y$, $Z(1) $  is distributed according to $\mathcal{M}\big( N; s(y) \big) $  a Multinomial with infinitely many classes. Then conditionally to $\{ \tau = 1 \}$, the probability that there is some particle at $-1$ is given by
$$
\frac{s_{1}(y) }{(1-s_0(y))} \, .
$$
The positions of the particles remain independent under the conditional probability and we conclude that
$$
\P_y \big(Z_1(1) = \cdot \mid \tau=1  \big) = \mathcal{B}\big(N; s_{1}(y)/(1-s_0(y)) \big)(\cdot) \,.
$$
Assuming that $y_0 \geq 1$ otherwise we must consider $\tilde{y}$ the shifted vector
\begin{align*}
\frac{s_{1}(y)}{1-s_0(y)}  & \geq \frac{\Bigr( 1- p_0 \Bigl)^{y_0}  - q_2^{y_{0}} }{\Bigr( 1- p_0 \Bigl)^{y_0}}\\
& \geq 1 - (\theta'\,)^{y_{0}} > \alpha ,
\end{align*}
where the lower bound holds for $N$ large enough as consequence of Assumption (A') and the definition of $\alpha$ . A large deviation argument allow us to conclude that for $\varepsilon$ small enough 
$$
\P_y \big( Z_{1}(1) \in \Xi \mid Z_0(1) = 0  \big) \geq \P_y \big( Z_{1}(1) \geq ( \alpha + \varepsilon ) N \mid Z_0(1) = 0  \big) \to 1 \,.
$$ 
Then, the infimum in (\ref{inq:4}) is larger than any $c_{\arabic{aux}}<1$ for $N$ sufficiently large. This finishes the proof. \hfill $\Box$
\vspace{.2cm}

The next Corollary generalizes (\ref{equa.aim.section.bound}) to the latter visiting times of $\triangle$.

\begin{cor} \label{cor:taui,bouded} Assume that $\xi$ satisfies Assumption (A'). Then, for every $i \in \N$, $\sup_{x \in \Omega }\E_x [T_\triangle^{(i)} ]$ and $\sup_{x \in \Omega }\E_x [\tau^{(i)} ]$ are bounded uniformly on $N$. In particular, under $\P_x$ the families of random variables $T_\triangle^{(i)}$ and $\tau^{(i)}$ are uniformly integrable.
\end{cor}
\emph{Proof.} Since $\tau^{(i)} \leq T_\triangle^{(i)}$, it suffices to prove the statements for $T_\triangle^{(i)}\,.$ To prove that the expectation is bounded we proceed inductively and apply the strong Markov property at time $T_\triangle^{(i-1)}\,.$ It is clear that
$$
\sup_{x \in \Omega }\E_x [T_\triangle^{(i)} ] \leq K^i \,,
$$ 
where $K$ is given by (\ref{equa.prop:bound.integral.T}). We now prove the uniform integrability. Applying the Markov property we obtain the upper bound
$$
\E_x [ T_\triangle^{(i)} \,;\, T_\triangle^{(i)} \geq l ] \leq \big(  \sup_{x\in \Omega(N)}  \E_{x} [ T_\triangle^{(i)} ] +l \, \big) \P_x  ( T_\triangle^{(i)} \geq l  )\,.
$$
It is not difficult to see that the right-hand side of the Equation converges to zero, finishing the proof.
\hfill $\Box$




\subsection{Convergence of Some Related Integrals}\label{sec.conv.integral}

To compute the front velocity in Subsection \ref{sec:front.speed.3.states}, we have to calculate two integrals $ \E_\oplus [T_\triangle]$ and $\E \big[ \phi \big( X ( T_\triangle ) \big) \big]$, where in the latter we assume that all particles start from zero. Hence
$$
\phi \big( X ( T_\triangle ) \big) = -  \sum_{i=1}^{\infty} \min \{ l \in \N ; Z_l (\tau^{(i)}) \neq 0 \} 1_{\{\tau^{(i)} \leq T_\triangle \}} \,.
$$ 
In the next Lemma, we use for the first time the condition $\E[|\vartheta|] < \infty$ that appears in Assumption (A) and Assumption (A').  
  
\begin{lem}\label{lem.int.conv.T.triangle.and.phi.X.eta>0}  Assume that $\xi$ satisfies Assumption (A'). Then for every $x \in \Omega(N)$
\begin{equation}\label{equa.lem.int.conv.T.triangle.and.phi.X.eta>0}
 \E_x \big[ \min \{ l \in \N ; Z_l (\tau) \neq 0 \} \big] = 1+o(1)\,.
\end{equation}
The term $o(1)$ converges to zero, as $N \to \infty$, independently from the initial condition $x \in \Omega(N)$.
\end{lem}
\emph{Proof.} To prove the Lemma, it suffices to show that the left-hand side in (\ref{equa.lem.int.conv.T.triangle.and.phi.X.eta>0}) is bounded from above by $1+o(1)$. By an argument similar to the one used in Lemma \ref{lem:jumping.Lambda} we obtain that
\begin{align}\label{equa.prova.lem.int.conv.T.triangle.and.phi.X.eta>0}
\E_x &\Big[ \min\{ l \in \N; Z_l (\tau) \neq 0 \} \Big] \nonumber \\
& \leq \sup_{y \in \Omega(N)}  \E_y \Big[ \min  \{l \in \N; Z_l (1) \neq 0 \} \mid \tau = 1 \Big] \nonumber \\
& =1 + \sup_{y \in \Omega(N)} \E_y \Big[ \min \{l \in \N;  Z_l (1) \neq 0 \} 1_{ \{ \min\{l \in \N; Z_l (1) \neq 0 \} \geq 2 \} } \mid \tau = 1 \Big]\,.
\end{align}
Under the conditional probability $Z$ is a Multinomial with infinitely many classes and parameters $s_i(y)/ \big(1-s_0(y) \big).$ Therefore, the probability that there is some particle at $-1$ is larger than $1-\theta'$, as a consequence of Assumption (A'). Moreover, the minimum is bounded from above by some $|\xi_{ij}|$. Indeed it suffices to choose $i$ such that $X_i(0)$ is a leader. Then
$$
-\min \{l \in \N;  Z_l (1) \neq 0 \} = \phi\big( X (1)\big) - X_i(0) \geq  \xi_{ij}.
$$
Hence, we can give an upper bound for the right-hand side in (\ref{equa.prova.lem.int.conv.T.triangle.and.phi.X.eta>0}). In fact, for $y \in \Omega(N)$ we obtain that
\begin{align*}
\E_y & \Big[ \min \{ l \in \N; Z_l (1) \neq 0 \} 1_{ \{ \min \{l \in \N; Z_l (1) \neq 0 \} \geq 2 \} } \mid \tau = 1 \Big] \\
&\leq \E_y \Big[| \xi_{ij} | 1_{ \{ \min \{ l \in \N; Z_l (1) \neq 0 \} \geq 2 \} } \mid \tau = 1 \Big] \\
& =  \E\big[ |\vartheta_{ij}| \big] \P_y \Big( \min \{ l \in \N; Z_l (1) \neq 0 \} \geq 2  \mid \tau= 1 \Big) \leq \E\big[ |\vartheta_{ij}| \big] (\theta')^N\,.
\end{align*}
It converges to zero independently from the initial position $y \in \Omega(N)$. \hfill $\Box$

With Lemma \ref{lem.int.conv.T.triangle.and.phi.X.eta>0} at hand we prove the following result in the noncritical case.

\begin{prop}\label{prop.int.conv.T.triangle.and.phi.X.eta>0}  Assume that $\xi$ satisfies Assumption (A') with $\eta>0\,$. Then 
\begin{equation}
\lim_{N\to \infty} \E_\oplus [T_\triangle ] = (m+1)  \, \quad and \quad \lim_{N \to \infty }\E \big[ \phi \big(\,X ( T_\triangle )  \big) \big] = -1\,.
\end{equation}
\end{prop}
\emph{Proof.} The first limit is a direct consequence of the uniform integrability of $T_\triangle$ and $ \P_\oplus ( T_\triangle = m+1 ) \to 1 \,,$ as $N\to \infty\,$. We now prove the second statement.
\begin{align} \label{equa.prop.int.conv.T.triangle.and.phi.X.eta>0}
\E &\big[ \phi \big(  X( T_\triangle  )  \big)  \big] \nonumber \\
& = -\sum_{i=1}^{\infty} \E_\oplus \bigg[  \min \{l \in \N; Z_l (\tau^{(i)}) \neq 0 \} 1_{\{ T_\triangle \geq \tau^{(i)} \}}  \bigg]  \nonumber \\
& = -\sum_{ i =1 }^{\infty} \sum_{y \in \Omega(N)} \E_\oplus \bigg[  \E_y \Big[  \min \{ l \in \N; Z_l (\tau^{(i)}) \neq 0 \} \Big] 1_{\{Z(\tau^{(i-1)})= y \,;\, T_\triangle \geq \tau^{(i)} \}} \bigg]   \nonumber \\
& = \big( -1 +o(1) \big) \sum_{i=1}^{\infty} \P_\oplus (T_\triangle \geq \tau^{(i)}) \, .
\end{align}
The last equality in (\ref{equa.prop.int.conv.T.triangle.and.phi.X.eta>0}) is a consequence of Lemma \ref{lem.int.conv.T.triangle.and.phi.X.eta>0}. The sum in (\ref{equa.prop.int.conv.T.triangle.and.phi.X.eta>0}) also converges to one. Indeed, it is a consequence of the strong Markov property and the uniform integrability of $T_\triangle\,.$ Hence, we obtain that the Equation in (\ref{equa.prop.int.conv.T.triangle.and.phi.X.eta>0}) converges to one, which finishes the proof of the Proposition. \hfill $\Box$

The critical case is more delicate and we prove the following result.

\begin{prop}\label{prop.int.conv.T.triangle.and.phi.X}  Assume that $\xi$ satisfies Assumption (A) and that $\eta=0\,$. Then 
\begin{equation}\label{equa.prop.int.conv.T.triangle.and.phi.X}
\lim_{N\to \infty} \E_\oplus [T_\triangle ] = (m+1)\E_0 [\mathcal{T}_0 ] - 1  \quad and \quad \lim_{N \to \infty }\E [\phi \big(X ( T_\triangle ) \big)] = - \E_0 [\mathcal{T}_0]\,,
\end{equation}
where $\mathcal{T}_0$ is the stopping time given by $\mathcal{T}_0 := \min\{ i \in \N; V^{(i)} = 0\}$, for $V^{(i)}$ a Markov chain defined as in Proposition \ref{prop:conv.joint.law.process}.
\end{prop}
\vspace{.2cm}

From Proposition \ref{prop:conv.joint.law.process}, $Z_1(\tau^{(i)})/N  $ converges in distribution to $G(V^{(i)}) $ as $N$ goes to infinity. We would like to state that
$$
T_\triangle = \min \{ i \in \N; Z_1(\tau^{(i)})/N = 1 \}  \stackrel{d}{\rightarrow} \min \{ i \in \N; G(V^{(i)}) = 1 \} \,. 
$$
Nevertheless the functional is not continuous and the above convergence must be justified. 

\begin{lem}\label{lem.conv.mathfrak.T} Assume that $\xi$ satisfies Assumption (A) with $\eta=0\,$. Then
$$
\min \{ i \in \N ; Z(\tau^{(i-1)} +m) = 0\} \stackrel{d}{\longrightarrow} \min \{ i \in \N ; V^{(i)} = 0 \} =\mathcal{T}_0 \,.
$$

\end{lem} 
\emph{Proof.} The minimum becomes continuous when restricted to $\N^\N\,.$ Since $ Z(\tau^{(i-1)} +m)$ converges in distribution to $V^{(i)}$ we conclude that the minimums also converge in distribution, which proves the result. We refer to \cite{Billingsley2009} for more details on convergence in distribution. \hfill $\Box$
\vspace{.2cm}

\begin{lem} \label{lem:law.convergence.T.triangle} Assume that $\xi$ satisfies Assumption (A) with $\eta=0\,$. Then $ T_\triangle $ converges in distribution to $\mathcal{T}_0 $ as $N \to \infty$. In particular,
$$
 T_\triangle \stackrel{d}{\to} \min \{ i \in \N ; G(V^{(i)}) = 1  \} \,.
$$
\end{lem}
\emph{Proof.} From Proposition \ref{prop:chain.started.random.position} 
$$
\min \{ i \in \N ; Z_0(\tau^{(i-1)} +m) = 0\}  -  \min \{ i \in \N ; Z_1(\tau^{(i)})/N  =1 \} \,,
$$ 
converges in distribution to zero. Hence, from Lemma \ref{lem.conv.mathfrak.T} we obtain that $T_\triangle$ converges in distribution to $\mathcal{T}_0 $. The second statement follows from
$$
\min \{ i \in \N ; G(V^{(i)}) =1  \} = \min \{ i \in \N ; V^{(i)}= 0 \} = \mathcal{T}_0 \quad a.s.
$$ \hfill $\Box$
\vspace{.2cm}


\emph{Proof of Proposition \ref{prop.int.conv.T.triangle.and.phi.X}.} It is not hard to see that if $\min \{ i \in \N ; Z_1(\tau^{(i)})/N =1 \} = k $, then $T_\triangle = \tau^{(k)}$. So we can write
\begin{align*}
 \E_\oplus [T_\triangle ] & = \sum_{k=1}^{\infty} \E_\oplus [  \tau^{(k)} \, 1_{\{T_\triangle = \tau^{(k)} \}} ] \\
 & = \sum_{k=1}^{\infty} \E_\oplus \bigg[ \sum_{j=1}^k \big(\tau^{(j)} - \tau^{(j-1)} \big)\, 1_{\{ \min_{i \in \N } \{Z_1(\tau^{(i)})/N  =1 \}  = k \}}  \bigg] \, .
\end{align*}
For a fixed $k$ the random variable $ \sum_{j=1}^k ( \tau^{(j)} - \tau^{(j-1)}  ) \,1_{\{ \min_{i \in \N } \{  Z_1(\tau^{(i)})/N = 1\}  = k \}}$ converges in law to 
$$
\sum_{j=1}^k (m + 1_{\{ V^{(j)} \neq 0 \}}\, )\, 1_{\{ \min \{ i \in \N ;   G(V^{(i)}) =1  \}  = k \}} = \big(( m+1 ) \, k  - 1 \big) 1_{\{ \mathcal{T}_0 = k \}} \,.
$$ 
Since $\tau^{(k)}$ is uniformly integrable the convergence also holds in $L^1$. From the uniform integrability of $T_\triangle$ we obtain the convergence in $L^1$ of the sum and the following limit holds. 
\begin{align*}
\lim_{N\to \infty} \E_\oplus [T_\triangle ] & = \sum_{k=1}^{\infty} \big( (m+1) k -1 \big) \P_0 \big(  \mathcal{T}_0 = k  \big) \notag  \\
 & =  (m+1) \E_0 \big[  \mathcal{T}_0  \big]  -\sum_{k=1}^{\infty} \P_0 \big( \mathcal{T}_0 = k \big)  \\
&= (m+1) \E_0 \big[  \mathcal{T}_0 \big]  - 1\, .
\end{align*}
This proves the first statement of Proposition \ref{prop.int.conv.T.triangle.and.phi.X}. We now prove the second limit in (\ref{equa.prop.int.conv.T.triangle.and.phi.X}). From the proof of Proposition \ref{prop.int.conv.T.triangle.and.phi.X.eta>0} we obtain that
$$
\E_\oplus \big[ \phi \big(  X( T_\triangle  )  \big)  \big] = -\big(1 +o(1) \big) \sum_{i=1}^{\infty} \P_\oplus \big( \tau^{(i)} \leq T_\triangle \big)  
$$  
From the uniform integrability of $T_\triangle$ we obtain that $\sum_{i=1}^{\infty} \P_\oplus \big( \tau^{(i)} \leq T_\triangle \big) \to  \E_0 \big[  \mathcal{T}_0 \big]$ which finishes the proof.  \hfill $\Box$

\vspace{.2cm}

The transition matrix of $V^{(i)}$ depends on $G$ and a fortiori on $\theta$. A coupling argument shows that $E_0 [\mathcal{T}_0]$ is non-increasing in $\theta$. Nevertheless, we do not know how to calculate explicitly the integral. However the asymptotic behaviors as $\theta \to 0 \,$ and $1$ are easy to compute.

\begin{prop}\label{prop:conv.gamma.theta.to.0} Let $V^{(i)}$ be the Markov chain whose transition matrix is given in Proposition \ref{prop:conv.joint.law.process} , then
$$
\lim_{\theta \to 0}   E_0 [\mathcal{T}_0 ] = \exp \rho^m \, .
$$
\end{prop}
\emph{Proof.} We write $E_0 [ \mathcal{T}_0]= \sum_{k=0}^{\infty} P_0  ( \mathcal{T}_0 \geq k+1 ) \, .$ For $l \geq 1$, then $1 \geq G(l) \geq G(1) = 1-\theta$, and 
\begin{align*}
P_0 & (\mathcal{T}_0 \geq k+1 ) \\
& = \sum_{ l_1 = 1}^{\infty} e^{-\rho^{m}} \frac{(\rho^{m} )^{l_1}}{l_1 !} \ldots \sum_{ l_{k-1} = 1}^{\infty} e^{-\rho^{m} G(l_{k-2})} \frac{(\rho^{m} G(l_{k-2}) \,)^{l_{k-1}}}{l_{k-1} !} \big(1- e^{-\rho^{m} G(l_{k-1})} \big) \,.
\end{align*}
The last expression is bounded from above by $( 1- e^{-\rho^{m} } )^k$. Since $ G(l) \to 1 $ as $\theta \to 0 $ we can conclude by dominated convergence. \hfill $\Box$

\vspace{.2cm}

We point out here that the case $\theta \to 0$ corresponds to the two state model studied in Section \ref{sec.derrida.two}. Informally, when $\theta$ is very small there is a high probability that $Z(\tau)$ starts afresh from $\triangle$. A similar computation can be done in the case where $\theta$ converges to one

\begin{prop} Let $V^{(i)}$ be the Markov chain whose transition matrix is given in Proposition \ref{prop:conv.joint.law.process}. Then
$$
\lim_{\theta \to 1}  E_0 [\mathcal{T}_0 ] = 2 - \exp - \rho^m \, .
$$
\end{prop}

\emph{Proof.} The proof follows the same lines as that of Proposition \ref{prop:conv.gamma.theta.to.0} and we leave the details to the reader. \hfill $\Box$

\subsection{Front speed}\label{sec:front.speed.3.states}

As in Subsection \ref{subsec.front.speed}, we explore the renewal structure of $Z$ that starts afresh from $\triangle$. Let $N (t) = \max \{ i \ ; \ T_\triangle^{(i)} \leq t \} $. Then
$$
\phi \big( X(t)\big) =  - \sum_{i=1}^{N(t)} \Big[ \phi \big( X(T_\triangle^{(i+1)} ) \big) -\phi \big( X( T_\triangle^{(i)}) \big) \Big] + o (t) \, .
$$
Taking the limit, as $t \to \infty$, we have that
\begin{align} \label{equa:eta0.aim}
\lim_{t \to \infty }\frac{\phi \big( X(t) \big)}{t} & = \lim_{t \to \infty} - \frac{1}{t}\sum_{i=1}^{N (t)} \phi \big(  X( T_\triangle^{(i+1)} ) \big) - \phi \big( X( T_\triangle^{(i)})  \big) \nonumber \\
& =  \frac{\E \big[ \phi \big( X(T_\triangle ) \big) \big]}{\E_\oplus [T_\triangle ] } \qquad a.s.
\end{align}
The limit is a consequence of the ergodic Theorem and the renewal structure. In Subsection \ref{sec.conv.integral}, we calculated the limits of the above expected values. We obtain that
$$
\lim_{N \to \infty} v_N = \left\{ 
\begin{array}{lcl}
-\big( 1+ \lfloor  1/r  \rfloor  \big)^{-1}, & \mbox{if} & 1/r \not\in \N   \\
- \big( \lfloor   1/r   \rfloor  + 1- 1/E_0 [ \mathcal{T}_0] \big)^{-1} , & \mbox{if} & 1/r = m \in \N \,,
\end{array}
\right.
$$
which proves Theorem \ref{teo.speed.3.states} with $g(\theta) = E_0 [  \mathcal{T}_0]$.
\section{Conclusion and sketch of the proof of Theorem \ref{teo.speed.3.states.gen}} \label{sec.conclusion.teo.speed.3.states.gen}

Theorem \ref{teo.speed.3.states.gen} follows as a Corollary of Theorem \ref{teo.speed.3.states} proved in Section \ref{sec.front.speed.three.states}. We will not prove it in detail but we give a sketch of the proof. The constants $\lambda_0$ and $\lambda_1 - \lambda_0$ appearing in Theorem \ref{teo.speed.3.states.gen} are justified by an affine transformation. Then, it remains to explain how we pass from the distribution over the lattice to the more general one. In the proof of Theorem \ref{teo.speed.3.states} we see that in the discrete case $\vartheta $ contributes to the position of the leaders only in rare events. Indeed, if there are $k$ leaders at time $t$ the position of the front is determined by $\vartheta$ at $t+1$ only in the case where $\xi_{ij} (t+1) \leq -2 $ for at least $N^k$ random variables. The probability of this event is of order $\theta^{N^{k}}$, as a consequence of Assumption (A). This behavior still holds in the general case. For a complete proof we refer to \cite{Comets2013} Theorem 1.3, which applies also to our case with the obvious changes.
\medskip

The position of the front depends basically on the tail distribution of $\xi$, that is determined by the point masses $\lambda_0$ and $\lambda_1$. The only case where $\vartheta $ could contribute to the position of the front in long time scales is in the non-integrable case. Then the mechanism responsible for propagation is of a very different nature and the front is no longer pulled by the leading edge. In the rare events, when the front moves backwards more than $\lambda_0 - \lambda_1$ the contribution of $\vartheta$ would be non-negligible depending on its tail and the global front profile. This problem is still open and much harder to solve. 

\section*{Acknowledgments}
\emph{I thank Francis Comets for suggesting  this problem to me and for his guidance in my Ph.D.}

\nocite{*}
\bibliographystyle{plain}
\bibliography{front_velocity}
\vspace{1.0cm}
Universit\'e Paris Diderot - Paris 7, Math\'ematiques, case 7012, F-75 205 Paris Cedex 13, France.
\\
\emph{E-mail address:} \texttt{cortines@math.univ-paris-diderot.fr}

\end{document}